\documentclass{amsart}
\usepackage{amsmath,amsthm,amssymb,IMjournal}
\usepackage{graphicx}

\newcommand{\cf}{\emph{cf}}

\newcommand{\Real}{\mathbb{R}}

\newcommand{\Sphere}{\mathbb{S}}
\newcommand{\eps}{\varepsilon}
\newcommand{\sii}{L^2}
\newcommand{\Dom}{\mathfrak{D}}
\begin{document}
\newtheorem{Theorem}{Theorem}[section]

\numberwithin{equation}{section}

\title{Spectrum of the Laplacian in narrow tubular neighbourhoods
of hypersurfaces with combined Dirichlet and Neumann boundary conditions}

\author{\|David |Krej\v{c}i\v{r}\'{\i}k|, \v{R}e\v{z}}

\abstract 
We consider the Laplacian in a domain squeezed
between two parallel hypersurfaces in Euclidean spaces of any dimension,
subject to Dirichlet boundary conditions on one of the hypersurfaces
and Neumann boundary conditions on the other.
We derive two-term asymptotics for eigenvalues
in the limit when the distance between the hypersurfaces tends to zero.
The asymptotics are uniform and local in the sense that
the coefficients depend only on the extremal points where
the ratio of the area of the Neumann boundary to the Dirichlet one 
is locally the biggest.
\endabstract
\thanks
The research has been supported by 
RVO61389005 and the GACR grant No.\ P203/11/0701.
\newline
Written for proceedings of Equadiff 2013, to appear in Mathematica Bohemica.
\emph{Date:} 11 Nov 2013 
\endthanks
%

\section{Introduction}
%
Let~$\Sigma$ be a connected orientable $C^2$ hypersurface 
(compact or non-compact) in~$\Real^d$, with $d \geq 2$,
equipped with the Riemannian metric~$g$ induced by the embedding. 
The orientation is specified by a globally defined 
unit normal vector field $n:\Sigma\to \Sphere^{d-1}$.
Given a small positive parameter~$\eps$,
we consider the tubular neighbourhood
\begin{equation}\label{layer.intro}
  \Omega_\eps := \big\{x+\eps\,t\,n(x) \in \Real^d \ \big| \ 
  (x,t) \in \Sigma \times (0,1) \big\}
  \,.
\end{equation}
We always assume that the map $(x,t) \mapsto x+\eps\,t\,n(x)$
is injective on $\overline{\Sigma} \times [0,1]$;
in particular, we require that the principal curvatures of~$\Sigma$,
$\kappa_1,\dots,\kappa_{d-1}$, are bounded functions.
Let $-\Delta_\textit{DN}^{\Omega_\eps}$ be the Laplacian on~$\Omega_\eps$,
subject to Dirichlet and Neumann boundary conditions
on~$\Sigma$ and $\Sigma_\eps:=\Sigma+\eps\,n(\Sigma)$, respectively. 
If the boundary~$\partial\Sigma$ is not empty,
we impose Dirichlet boundary conditions on the remaining
part of~$\partial\Omega_\eps$.
We arrange the eigenvalues below the essential spectrum
of $-\Delta_\textit{DN}^{\Omega_\eps}$
in an increasing order and repeat them according to multiplicity,
$
  \lambda_1(\eps) \leq \lambda_2(\eps) \leq \lambda_3(\eps) \leq \dots 
$,
with the convention that all eigenvalues are included
if the essential spectrum is empty.
In fact, we make the sequence always infinite by defining 
$\lambda_n := \inf\sigma_\mathrm{ess}(-\Delta_\textit{DN}^{\Omega_\eps})$
for all $n>N$, if the number of eigenvalues below the essential spectrum
is a finite (possibly zero) natural number~$N$. 

\begin{figure}[h!]
\begin{center}
\includegraphics[width=\textwidth]{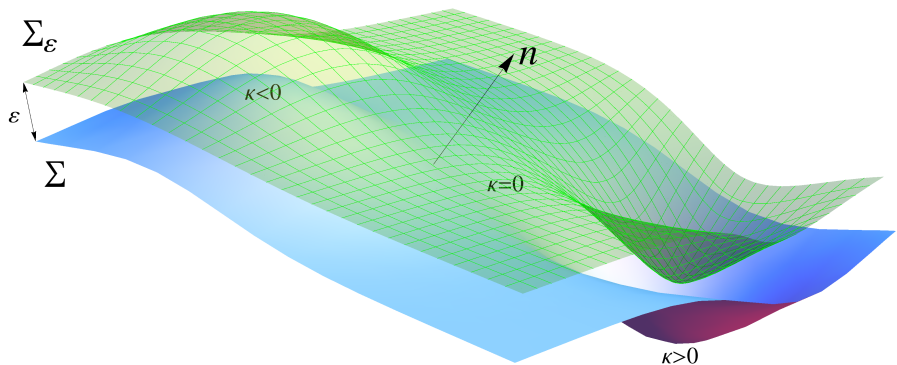}
\vspace{-25ex}
\end{center}
\caption{The geometry of the tubular neighbourhood~$\Omega_\eps$ 
for $d=3$.}\label{Fig}
\end{figure}

The objective of this paper is to show that 
the $d$-dimensional differential operator $-\Delta_\textit{DN}^{\Omega_\eps}$
can be approximated in the limit as $\eps \to 0$ 
by the $(d-1)$-dimensional Schr\"odinger-type operator
\begin{equation}\label{op.comparison}
  H_\eps := -\Delta_g  + \frac{\kappa}{\eps} 
  \qquad \mbox{on} \qquad
  \sii(\Sigma)
  \,.
\end{equation}
Here~$-\Delta_g$ denotes the Laplace-Beltrami operator of~$\Sigma$,
subject to Dirichlet boundary conditions if~$\partial\Sigma$ is not empty, 
and 
$
  \kappa := \kappa_1+\dots+\kappa_{d-1}
$ 
is a $d-1$ multiple of the mean curvature of~$\Sigma$.
Note that the sign of~$\kappa$ depends on the choice of orientation~$n$, 
that is on the direction in which the parallel surface~$\Sigma_\eps$ 
is constructed with respect to~$\Sigma$, 
\cf~Figure~\ref{Fig}.
We arrange the eigenvalues below the essential spectrum
of the operator~$H_\eps$ using the same conventions as above,
$
  \mu_1(\eps) \leq \mu_2(\eps) \leq \mu_3(\eps) \leq \dots 
$.

In this paper we establish the following spectral asymptotics:
\begin{Theorem}\label{Thm.main}
For all $n \geq 1$,
\begin{equation}\label{expansion}
  \lambda_n(\eps) 
  = \left(\frac{\pi}{2\eps}\right)^2 + \mu_n(\eps) + \mathcal{O}(1)
  \qquad \mbox{as} \qquad
  \eps \to 0
  \,.
\end{equation}
\end{Theorem}
%


This asymptotic expansion was proved previously
by the author for~$d=2$ in~\cite{K5}.
Moreover, some form of norm-resolvent convergence of 
$-\Delta_\textit{DN}^{\Omega_\eps}$ to~$H_\eps$ was established
and the result~\eqref{expansion} for $d=3$ was announced there.
In the present paper we extend the validity of formula~\eqref{expansion}
to any dimension and provide some details of the variational  proof 
which were missing in~\cite{K5}.

Using known results about the strong-coupling/semiclassical asymptotics
of eigenvalues of the Schr\"odinger-type operator~\eqref{op.comparison},
one has, for all $n \geq 1$,
\begin{equation}\label{strong}
  \mu_n(\eps)
  = \frac{\inf\kappa}{\eps} + o(\eps^{-1})
  \qquad\mbox{as}\qquad
  \eps \to 0
  \,.
\end{equation}
This result seems to be well known;
we refer to~\cite[App.~A]{FK1} for a proof in a general Euclidean case,
which extends to the present situation.

Combining~\eqref{expansion} with~\eqref{strong},
we see that the two leading terms 
in the $\eps$-expansion of~$\lambda_n(\eps)$ are independent of~$n$. 
Furthermore, the geometry of~$\Omega_\eps$ is seen 
in these terms only \emph{locally},
through the minimal value of the mean curvature of~$\Sigma$.
In view of the leading role of the mean curvature~$\kappa$
in the surface element of~$\Sigma_\eps$, \cf~\eqref{h.fomula},
we see that the minimal values of the mean curvature on~$\Sigma$
corresponds to points for which, roughly, the Neumann boundary has 
``locally the largest area'' with respect to the opposite Dirichlet one;
see also~Figure~\ref{Fig}. 
The results \eqref{expansion}--\eqref{strong} are thus consistent 
with the physical intuition that ``Dirichlet conditions raise energies
and Neumann conditions lower energies''.

The particular form of the thin-width expansions~\eqref{expansion} 
has important 
physical consequences for spectral properties of quantum waveguides
as explained in~\cite{K5}.
Let us also mention that the local character resembles 
situations of Dirichlet tubes of variable radius
\cite{Friedlander-Solomyak_2007,Friedlander-Solomyak_2008a,
Borisov-Freitas_2009,Borisov-Freitas_2010,Lampart-Teufel-Wachsmuth}.

The case of Neumann or Dirichlet tubes of uniform radius
differs from the present situation in many respects.
Let us denote by $\{\lambda_n^N(\eps)\}_{n=1}^\infty$
and $\{\lambda_n^D(\eps)\}_{n=1}^\infty$
the set of eigenvalues below the essential spectrum
of the Neumann and Dirichlet Laplacian
on $\sii(\Omega_\eps)$, respectively,
with the same conventions as used above for $\{\lambda_n(\eps)\}_{n=1}^\infty$.
The case of the Neumann Laplacian is trivial in the sense that
its spectrum is known to converge to the spectrum 
of the the underlying manifold~$\Sigma$, \cf~\cite{Schatzman_1996}.
More precisely,
\begin{equation}\label{expansion.Neumann}
  \lambda_n^N(\eps) 
  = 0 + \mu_n^N + o(1)
  \qquad \mbox{as} \qquad
  \eps \to 0
  \,,
\end{equation}
where $\{\mu_n^N\}_{n=1}^\infty$ is the set of eigenvalues 
below the essential spectrum (with the aforementioned conventions)
of the Laplace-Beltrami operator~$-\Delta_g$ on~$\sii(\Sigma)$,
subject to Neumann boundary conditions on~$\partial\Sigma$.
In order to consistently compare~\eqref{expansion.Neumann}
with~\eqref{expansion} (and~~\eqref{expansion.Dirichlet} below),
we included into~\eqref{expansion.Neumann}
the vanishing lowest Neumann eigenvalue 
of the transverse interval~$(0,\eps)$ 
and will refer to~$\mu_n^N$ as the ``second term'' 
in the expansion of~$\lambda_n^N(\eps)$.
In the Dirichlet case, we have~\cite{KRT}
\begin{equation}\label{expansion.Dirichlet}
  \lambda_n^D(\eps) 
  = \left(\frac{\pi}{\eps}\right)^2 + \mu_n^D + \mathcal{O}(1)
  \qquad \mbox{as} \qquad
  \eps \to 0
  \,,
\end{equation}
where $\{\mu_n^D\}_{n=1}^\infty$ is the set of eigenvalues 
below the essential spectrum 
(again with the aforementioned conventions)
of the Schr\"odinger-type operator 
$-\Delta_g  + V_\mathrm{eff}$ on~$\sii(\Sigma)$,
subject to Dirichlet boundary conditions on~$\partial\Sigma$.
Here~$V_\mathrm{eff}$ is a purely geometric, $\eps$-independent potential,
expressed solely in terms of the principal curvatures,
\begin{equation}\label{V.eff}
  V_{\mathrm{eff}}
  :=-\frac{\kappa_1^2+\dots+\kappa_{d-1}^2}{2} 
  +\frac{(\kappa_1+\dots+\kappa_{d-1})^2}{4}
  \,.
\end{equation}
Summing up, contrary to Theorem~\ref{Thm.main},
in the purely Neumann or Dirichlet case
the second term in the asymptotic expansion of eigenvalues
is independent of~$\eps$
and determined by the \emph{global} geometry of~$\Sigma$.

In addition to this introductory part, 
the paper consists of Section~\ref{Sec.Pre},
in which we collect some auxiliary material,
and Section~\ref{Sec.proof} devoted 
to the proof of Theorem~\ref{Thm.main}.

\section{Preliminaries}\label{Sec.Pre}
%
We refer to~\cite{KRT} for a necessary geometric background 
of tubes about hypersurfaces. 
Using the Fermi ``coordinates''~$(x,t)$ that appear in~\eqref{layer.intro},
$\Omega_\eps$~can be identified with the Riemannian manifold $\Sigma\times(0,1)$ 
equipped with the metric~$G$ of the following block-diagonal structure 
$G = G_{\mu\nu} \, dx^\mu dx^\nu + \eps^2 dt^2$.
Here the range of Greek indices is assumed to be $1,\dots,d-1$
and the Einstein summation convention is employed.
We shall not need the explicit formulae for the coefficients~$G_{\mu\nu}$,
just the bounds:
\begin{equation}\label{eq:metric_bound}
  (1-C\eps)(g_{\mu\nu})\leq(G_{\mu\nu})\leq (1+C\eps)(g_{\mu\nu})
  \,.
\end{equation}
(Of course, we implicitly assume that~$\eps$ is so small
that $1-C\eps$ is positive.)
Here and in the sequel, we adopt the convention that~$C, c$ 
and the constants involved in the ``big~$\mathcal{O}$'' notation
possibly depend on the supremum norm of the principal curvatures
$\kappa_1,\dots,\kappa_{d-1}$ and may vary from line to line.
On the other hand, we shall need the formula
for the determinant $|G| = \varepsilon^2 \, |g| \, h_\eps^2$,
where
\begin{equation}\label{h.fomula}
  h_\eps(\cdot,t) := \prod_{\mu=1}^{d-1}(1-\varepsilon \, \kappa_{\mu} \, t)
  = 1 - \eps \, \kappa \, t + \mathcal{O}(\eps^2)
  \,.
\end{equation}
The volume element of~$\big(\Sigma\times(0,1),G\big)$ 
is thus given by $d\Omega_\eps = \eps \, h_\eps \, d\Sigma \wedge dt$,
where $d\Sigma = |g|^{1/2} dx^1 \wedge \dots \wedge dx^{d-1}$
is the surface element of $(\Sigma,g)$.

Using the above geometric preliminaries,
the Hilbert space $\sii(\Omega_\eps)$ can be identified  with 
$
  \mathcal{H}_\eps :=
  \sii\big(\Sigma\times(0,1),\eps \, h_\eps \, d\Sigma \wedge dt\big)
$.
The Laplacian $-\Delta_\textit{DN}^{\Omega_\eps}$ 
can be in turn identified with the self-adjoint operator 
on $\mathcal{H}_\eps$ associated with the quadratic form 
\begin{align*}
  Q_\eps[\psi]
  &:= \big\langle
  \partial_{x^\mu}\psi,G^{\mu\nu}\partial_{x^\nu}\psi
  \big\rangle_{\mathcal{H}_\eps}
  + \eps^{-2} \|\partial_t\psi\|_{\mathcal{H}_\eps}^2
  \,,
  \\
  \psi \in \Dom(Q_\eps) &:=
  \left\{
  \psi \in W^{1,2}\big(\Sigma \times (0,1)\big) \ | \quad
  \psi = 0 \quad \mbox{on} \quad
  \partial\big(\Sigma \times (0,1)\big) \setminus \big(\Sigma\times\{1\}\big)
  \right\}
  \,.
\end{align*}
Here the boundary values of~$\psi$ are understood in the sense of traces.
Similarly, the operator~$H_\eps$ is associated with the form 
\begin{align*}
  q_\eps[\varphi]
  &:= \big\langle
  \partial_{x^\mu}\varphi,g^{\mu\nu}\partial_{x^\nu}\varphi
  \rangle_{\sii(\Sigma)}
  + \eps^{-1} \langle
  \varphi,\kappa\varphi
  \rangle_{\sii(\Sigma)}
  \,,
  \\
  \varphi \in \Dom(q_\eps) &:= W_0^{1,2}(\Sigma)
  \,.
\end{align*}

The spectral numbers $\{\lambda(\eps)\}_{n=1}^\infty$ as defined above
can be fully characterised by the Rayleigh-Ritz variational 
formula~\cite[Sec.~4.5]{Davies}
\begin{equation}\label{minimax} 
  \lambda_n(\eps) = 
  \inf_{\mathcal{L}_n} 
  \sup_{\psi \in \mathcal{L}_n}  
  \frac{Q_\eps[\psi]}{\ \|\psi\|_{\mathcal{H}_\eps}^2}
  \,,
\end{equation}
where the infimum is taken over all $n$-dimensional subspaces 
$\mathcal{L}_n\subset\Dom(Q_\eps)$.
An analogous formula holds for the spectral numbers
$\{\mu(\eps)\}_{n=1}^\infty$ of~$H_\eps$.
It follows from~\eqref{minimax} that the presence
of the multiplicative factor~$\eps$ in the weight of~$\mathcal{H}_\eps$
has no effect on the spectrum of~$-\Delta_\textit{DN}^{\Omega_\eps}$ .

Our strategy to prove Theorem~\ref{Thm.main} will be to show
that the forms~$Q_\eps$ and~$q_\eps$ are close to each other
in a sense as $\eps \to 0$. 
Since the forms act on different Hilbert spaces,
this requires a suitable identification of~$\mathcal{H}_\eps$ with~$\sii(\Sigma)$.
First, notice that it follows from~\eqref{h.fomula} that~$\mathcal{H}_\eps$ 
(up to the irrelevant factor~$\eps$)
approaches the $\eps$-independent Hilbert space 
$
  \mathfrak{H} :=
  \sii\big(\Sigma\times(0,1), d\Sigma \wedge dt\big)
$.
For this Hilbert space, we use the orthogonal-sum decomposition
\begin{equation}\label{direct}
  \mathfrak{H} = \mathfrak{H}_1 \oplus \mathfrak{H}_1^\bot
  \,,
\end{equation}
where the subspace~$\mathfrak{H}_1$ consists of functions~$\psi_1$
such that
\begin{equation}\label{psi1}
  \psi_1(x,t) = \varphi(x) \chi_1(t)
  \qquad \mbox{with} \qquad
  \varphi \in \sii(\Sigma)
  \,, \quad
  \chi_1(t):=\sqrt{2} \sin\left(\pi t/2\right)
  \,.
\end{equation}
Notice that~$\chi_1$ is the first eigenfunction of 
the Laplacian on $\sii((0,1))$, 
subject to the Dirichlet and Neumann boundary condition
at~$0$ and~$1$, respectively.
This operator has eigenvalues $\{(n\pi/2)^2\}_{n=1}^\infty$,
where the lowest one is of course related 
to the leading term in~\eqref{expansion}.
Since~$\chi_1$ is normalised, we clearly have
$
  \|\psi_1\|_{\mathfrak{H}}=\|\varphi\|_{\sii(\Sigma)}
$.
Given any $\psi\in\mathfrak{H}$, we have the decomposition
\begin{equation}\label{psi.decomposition}
  \psi = \psi_1 + \psi_\bot
  \qquad\mbox{with}\qquad
  \psi_1 \in \mathfrak{H}_1, \ \psi_\bot\in \mathfrak{H}_1^\bot
  \,,
\end{equation}
where~$\psi_1$ has the form~\eqref{psi1}
with $\varphi(x):=\int_0^1 \psi(x,t) \chi_1(t)  \;\! dt$.
Note that $\psi_1,\psi_\bot\in\Dom(Q_\eps)$ if $\psi\in\Dom(Q_\eps)$.
The inclusion $\psi_\bot\in\mathfrak{H}_1^\bot$ means that
\begin{equation}\label{orth.identity1}
  \int_0^1 \psi_\bot(x,t) \, \chi_1(x) \, dt = 0
  \qquad\mbox{for a.e.}\quad x \in \Sigma
  \,.
\end{equation}
If in addition $\psi_\bot \in \Dom(Q_\eps)$,
then one can differentiate the last identity to get
\begin{equation}\label{orth.identity2}
  \int_0^1 \partial_{x^\mu}\psi_\bot(x,t) \, \chi_1(t) \, dt = 0
  \qquad\mbox{for a.e.}\quad x \in \Sigma
  \,.
\end{equation}
Since $\mathcal{H}_\eps$ and~$\mathfrak{H}$ can be identified 
as vector spaces for any fixed~$\eps>0$,
the decomposition~\eqref{direct} 
can be equally used for each function~$\psi \in \mathcal{H}_\eps$. 
In view of the isomorphism 
$
  \sii(\Sigma)\ni\varphi \mapsto \psi_1 \in \mathfrak{H}_1
$,
we may think of~$H_\eps$ as acting on~$\mathfrak{H}_1$ as well.

\section{Proof of Theorem~\ref{Thm.main}}\label{Sec.proof}
%
Expansion~\eqref{expansion} will follow as a consequence
of upper and lower bounds to~$\lambda_n(\eps)$
that have the same leading order terms in their asymptotics.
It is convenient to define the shifted form 
$\tilde{Q}_\eps:=Q_\eps-\pi^2/(2\eps)^2$
and focus on the first non-trivial term~$\mu_n(\eps)$ in~\eqref{expansion}.
Let us decompose any $\psi \in \Dom(Q_\eps)$ 
according to~\eqref{psi.decomposition}.
A straightforward calculation employing an integration by parts yields
\begin{equation}\label{crucial}
\begin{aligned}
  \|\partial_t\psi\|_{\mathcal{H}_\eps}^2 
  - \left(\frac{\pi}{2}\right)^2 \|\psi\|_{\mathcal{H}_\eps}^2
  =\ & \|\partial_t\psi_\bot\|_{\mathcal{H}_\eps}^2 
  - \left(\frac{\pi}{2}\right)^2 \|\psi_\bot\|_{\mathcal{H}_\eps}^2
  - 2 \eps \, \Re \int \overline{\varphi}\, \chi_1' \,
  \psi_\bot \, \partial_t h_\eps 
  \\
  \ & + \frac{\eps}{2} \int |\varphi|^2 \, \chi_1^2 \, \partial_t^2 h_\eps 
  - \eps \int_\Sigma |\varphi|^2 \, \partial_t h_\eps |_{t=1} 
  \,.
\end{aligned}
\end{equation}
Here and in the sequel, 
$\int$ and $\int_\Sigma$
abbreviate the integrals over~$\Sigma\times(0,1)$ and $\Sigma$
with the integration measures~$d\Sigma \wedge dt$ and $d\Sigma$, respectively,
and we do not write the variables on which the integrated functions depend. 
Using~\eqref{h.fomula} and recalling that~$\chi_1$ is normalised, 
we easily verify
\begin{equation}\label{varphi.estimates}
\begin{aligned}
  \left|
  \frac{1}{\eps^2} \int |\varphi|^2 \, \chi_1^2 \, \partial_t^2 h_\eps 
  \right|
  \leq C \int_\Sigma |\varphi|^2
  \,, 
  \\
  \left|
  - \frac{1}{\eps^2} \int_\Sigma |\varphi|^2 \, \partial_t h_\eps |_{t=1} 
  - \eps^{-1} \big\langle
  \varphi,\kappa\varphi
  \rangle_{\sii(\Sigma)}
  \right|
  \leq C \int_\Sigma |\varphi|^2
  \,,
\end{aligned}
\end{equation}
which reveals the source of the potential term of~\eqref{op.comparison}. 
At the same time, using~\eqref{eq:metric_bound},
\begin{equation}\label{kinetic.estimates}
\begin{aligned}
  \pm \, \eps^{-1} \big\langle
  \partial_{x^\mu}\psi,G^{\mu\nu}\partial_{x^\nu}\psi
  \big\rangle_{\mathcal{H}_\eps}
  &\leq \pm (1 \pm C\eps) \,
  \big\langle
  \partial_{x^\mu}\psi,g^{\mu\nu}\partial_{x^\nu}\psi
  \big\rangle_{\mathfrak{H}}
  \,,
  \\
  \pm \, \eps^{-1} \|\psi\|_{\mathcal{H}_\eps}^2
  &\leq \pm (1 \pm C\eps) \, \|\psi\|_{\mathfrak{H}}^2
  \,.
\end{aligned}
\end{equation}
Here, by the normalisation of~$\chi_1$ 
and \eqref{orth.identity1}--\eqref{orth.identity2},
\begin{equation}\label{orth.identities}
\begin{aligned}
  \big\langle
  \partial_{x^\mu}\psi,g^{\mu\nu}\partial_{x^\nu}\psi
  \big\rangle_{\mathfrak{H}}
  &= \big\langle
  \partial_{x^\mu}\varphi,g^{\mu\nu}\partial_{x^\nu}\varphi
  \big\rangle_{\sii(\Sigma)}
  + \big\langle
  \partial_{x^\mu}\psi_\bot,g^{\mu\nu}\partial_{x^\nu}\psi_\bot
  \big\rangle_{\mathfrak{H}}
  \,,
  \\
  \|\psi\|_{\mathfrak{H}}^2
  &= \|\varphi\|_{\sii(\Sigma)}^2 + \|\psi_\bot\|_{\mathfrak{H}}^2
  \,.
\end{aligned}
\end{equation}

\subsection{Upper bound}
Let us restrict the subspaces~$\mathcal{L}_n$ in the formula~\eqref{minimax}
to the decoupled functions~\eqref{psi1}, where $\varphi \in \Dom(q_\eps)$. 
Using~\eqref{crucial}--\eqref{orth.identities} with $\psi_\bot=0$,
we get the upper bound
\begin{equation}\label{upper.pre}
  \frac{\tilde{Q}_\eps[\psi_1]}{\ \|\psi_1\|_{\mathcal{H}_\eps}^2}
  \leq \frac{(1+C\eps) \, q_\eps[\varphi]+C \, \|\varphi\|_{\sii(\Sigma)}^2}
  {(1-C\eps) \, \|\varphi\|_{\sii(\Sigma)}^2}
  \,,
\end{equation}
which yields
\begin{equation}\label{upper}
  \lambda_n(\eps) - \left(\frac{\pi}{2\eps}\right)^2
  \leq \frac{1+C\eps}{1-C\eps} \, \mu_n(\eps) + \frac{C}{1-C\eps}
  \,.
\end{equation}
Observing that, for each~$n \geq 1$,
\begin{equation}\label{mu.bounds}
  - \|\kappa\|_\infty \leq
  \eps \, \nu_n - \|\kappa\|_\infty
  \leq \eps \, \mu_n(\eps) \leq 
  \eps \, \nu_n + \|\kappa\|_\infty
  \,,
\end{equation}
where~$\nu_n$ are the ``eigenvalues'' of~$-\Delta_g$ as defined by~\eqref{minimax},
we conclude from~\eqref{upper} the desired asymptotic upper bound
\begin{equation}\label{upper.final}
  \lambda_n(\eps) 
  \leq \left(\frac{\pi}{2\eps}\right)^2 + \mu_n(\eps) + \mathcal{O}(1)
  \qquad \mbox{as} \qquad \eps \to 0 
  \,.
\end{equation}
It is worth noticing that the constant~$C$ in~\eqref{upper}
does not depend on~$n$; a possible dependence of the constants
appearing in~$\mathcal{O}(1)$ enters through 
the upper bound of~\eqref{mu.bounds} only.
 
\subsection{Lower bound}
As usual, lower bounds are more difficult to establish.
In our situation, we need to  carefully  
exploit the Hilbert-space decomposition~\eqref{direct}. 
Since $\psi_\bot \in \mathfrak{H}_1^\bot$, we have 
$
  \int_0^1 |\partial_t \psi_\bot(x,t)|^2 \, dt
  \geq \pi^2 \int_0^1 |\psi_\bot(x,t)|^2 \, dt
$
for a.e.\ $x \in \Sigma$. This Poincar\'e-type estimate 
extends to~$\mathfrak{H}$ by Fubini's theorem.
Hence, using~\eqref{eq:metric_bound} 
to estimate~$h_\eps$ in~$d\Omega_\eps$,
we get
\begin{equation*}
\begin{aligned}
  \|\partial_t\psi_\bot\|_{\mathcal{H}_\eps}^2 
  - \left(\frac{\pi}{2\eps}\right)^2 \|\psi_\bot\|_{\mathcal{H}_\eps}^2
  \geq 
  \eps \left[
  \left(\frac{3\pi^2}{4\eps^2}\right) 
  - C \left(\frac{5\pi^2}{4\eps}\right)
  \right]
  \|\psi_\bot\|_\mathfrak{H}^2
  \geq 
  \eps \, \frac{c}{\eps^2} \, \|\psi_\bot\|_\mathfrak{H}^2
  \,,
\end{aligned}
\end{equation*}
where the second inequality holds with a positive constant~$c$
for all sufficiently small~$\eps$. 
Using~\eqref{h.fomula} and the Young inequality,
the last term on the right hand side of~\eqref{crucial}
can be estimated as follows
\begin{equation*}
\begin{aligned}
  \frac{1}{\eps^2} \left|
  2 \, \Re \int \overline{\varphi}\, \chi_1' \,
  \psi_\bot \, \partial_t h_\eps 
  \right|
  & \leq \frac{C}{\eps} \,
  2 \int |\varphi \, \chi_1' \,
  \psi_\bot |
  \leq 
  \frac{C^2}{\delta} \|\varphi \, \chi_1'\|_\mathfrak{H}^2
  + \frac{\delta}{\eps^2} \|\psi_\bot\|_\mathfrak{H}^2
\end{aligned}
\end{equation*}
with any positive~$\delta$.
Here 
$
  \|\varphi \, \chi_1'\|_\mathfrak{H}
  = (\pi/2) \;\! \|\varphi\|_{\sii(\Sigma)}
$.
Choosing~$\delta$ sufficiently small
and using \eqref{crucial}--\eqref{orth.identities},
we thus get the lower bound
\begin{equation}\label{lower.pre}
  \frac{\tilde{Q}_\eps[\psi]}{\ \|\psi\|_{\mathcal{H}_\eps}^2}
  \geq \frac{(1-C\eps) \, q_\eps[\varphi]-C \, \|\varphi\|_{\sii(\Sigma)}^2
  +c \, \eps^{-2} \;\! \|\psi_\bot\|_\mathfrak{H}^2}
  {(1+C\eps) \, \big(\|\varphi\|_{\sii(\Sigma)}^2
  +\|\psi_\bot\|_{\mathfrak{H}}^2\big)}
  \,.
\end{equation}
Here the numerator is in fact the quadratic form of 
an operator direct sum
$ 
  T_\eps \oplus T_\eps^\bot
$
with respect to the decomposition~\eqref{direct},
where $T_\eps := (1-C\eps) \, H_\eps - C$
and $T_\eps^\bot := c \, \eps^{-2}$.
In view of~\eqref{mu.bounds}, the spectrum of~$T_\eps^\bot$
diverges faster as $\eps \to 0$ than that of~$T_\eps$.
This enables us to conclude from~\eqref{lower.pre}
with help of~\eqref{minimax} that for any $n \geq 1$ 
there exist $C,c$ such that for all $\eps \leq c$, we have
\begin{equation}\label{lower}
  \lambda_n(\eps) - \left(\frac{\pi}{2\eps}\right)^2
  \geq \frac{1-C\eps}{1+C\eps} \, \mu_n(\eps) - \frac{C}{1+C\eps}
  \,.
\end{equation}
Using~\eqref{mu.bounds}, we conclude from~\eqref{lower} 
the desired asymptotic lower bound  
\begin{equation}\label{lower.final}
  \lambda_n(\eps) 
  \geq \left(\frac{\pi}{2\eps}\right)^2 + \mu_n(\eps) + \mathcal{O}(1)
  \qquad \mbox{as} \qquad \eps \to 0 
  \,.
\end{equation}

Combining~\eqref{lower.final} with~\eqref{upper.final},
we complete the proof of Theorem~\ref{Thm.main}.

%
\bibliography{bib}
\bibliographystyle{amsplain}

\medskip

{\small
{\em Authors' addresses}:
{\em David Krej\v{c}i\v{r}\'{\i}k},
Nuclear Physics Institute ASCR, 
\v{R}e\v{z}, Czech Republic;
e-mail: \texttt{krejcirik@\allowbreak ujf.cas.cz}.
}

\end{document}